\documentclass[10pt]{article}
\usepackage{amsmath, amsthm}
\usepackage[T2A]{fontenc}
\usepackage[utf8]{inputenc}
\usepackage{enumerate}
\usepackage{amssymb}
\usepackage{bbold}
\usepackage{color}
\newtheorem{theorem}{Theorem}[section]
\newtheorem{proposition}[theorem]{Proposition}
\newtheorem{lemma}[theorem]{Lemma}
\newtheorem{example}[theorem]{Example}
\newtheorem{corollary}[theorem]{Corollary}
\newtheorem{definition}[theorem]{Definition}
\newtheorem{question}[theorem]{Question}

\begin{document}

\title{On automatic continuity of operators from ordered to topological vector spaces}
\date{\empty}
\maketitle
\author{\centering{E. Y. Emelyanov$^{1}$ and S. G. Gorokhova$^{2}$\\
\small $1$ Sobolev Institute of Mathematics, Novosibirsk, Russia\\
\small $2$ Southern Mathematical Institute of the Vladikavkaz Scientific Center of the Russian
Academy of Sciences, Vladikavkaz, Russia}\\

\abstract{We study continuity and boundedness of order-to-topology bounded and order-to-topology continuous operators
from ordered to topological vector spaces. Several results on automatic continuity
of operators from ordered Fr\'echet spaces to topological vector spaces are included. 
Levi and Lebesgue operators especially are investigated.}

\vspace{3mm}

\vspace{3mm}
}

\section{Introduction and preliminaries}

\hspace{4mm}
Automatic continuity of operators enjoying some additional properties plays significant role in analysis. 
For example, a nice theory exists on automatic continuity of algebra homomorphisms (see, e.g. \cite{D-1978,D-1984}).
Another case is the automatic continuity of positive operators from Banach lattices to normed lattices \cite{AB-2006,AT-2007,AN-2009}.
Apparently, importance of boundedness of operators from a Banach space to a normed space lies in the Uniform Boundedness Principle for families of these operators.
This paper is devoted to the study of the question how far one can go with other kinds of operators from a Banach space to a topological vector space.
This question has attracted an attention in several recent papers 
(see, for example, \cite{AEG-2022,E4-2025,E1-2026,E2-2026,EEG1-2025,EGS-2026,G-2024,G-2026,JAM-2021,KTA-2025,ZSC-2023}).

Throughout, vector spaces are real, and operators are linear. When suitable, we abbreviate a topological vector space 
(normed space, locally convex space, ordered vector space, ordered normed space, ordered Banach space) by \text{\rm TVS} 
(\text{\rm NS}, \text{\rm LCS}, \text{\rm OVS}, \text{\rm ONS}, \text{\rm OBS} respectively).
In what follows, symbol $x_\alpha\downarrow 0$ stands for a decreasing net in an OVS such that $\inf\limits_\alpha x_\alpha=0$;
$[a,b]=\{x\in X: a\le x\le b\}$ for an order interval in an OVS $X$; ${\cal L}(X,Y)$ for the space of operators between vector spaces $X$ and $Y$;
$B_X(a,r)=\{x\in X: p(x,a)\le r\}$ for the closed ball centered at $a\in X$ of radius $r\ge 0$ in a metric space $(X,p)$; and $B_X:=B_X(0,1)$, where $X$ is a NS.
It is proved in \cite[Theorem 2.1]{E2-2026} that each order-to-norm bounded family of operators from a Banach lattice to an NS is equi-continuous.
In \cite[Theorem 2.1]{EEG1-2025}, this result is extended to the OBS setting. In \cite{E1-2026}, the further
extension is obtained in the case of a TVS in the co-domain. In the present paper, we prove Theorem \ref{theorem 1} that generalizes 
\cite[Theorem 2.1]{E1-2026} to operators from an ordered Fr\'echet space with a closed generating cone to a TVS.

Let $T: X\to Y$ be an operator from a normed space $X$ to a TVS $(Y,\tau)$. 
In absence of additional structures in $X$, one cannot say much about conditions providing topological boundedness/continuity of $T$.
For example, using of finite sets gives no proper subclass, because $T$ always carries finite sets into $\tau$-bounded sets.
On the other hand, usage of compactness returns us to bounded operators. Indeed, if $T$ is $\tau$-bounded then $T$
takes compact sets into bounded sets. On the other hand, if $T$ is not $\tau$-bounded then there exist a circled $U\in\tau(0)$ 
and a sequence $(x_n)$ in $B_X$ such that $Tx_n\notin n^2U$ for all $n\in\mathbb{N}$. In this case, 
the image $\big\{T(\frac{x_n}{n})\big\}_{n=1}^\infty\cup\{0\}$ of the compact set $\big\{\frac{x_n}{n}\big\}_{n=1}^\infty\cup\{0\}$ 
under $T$ is not $\tau$-bounded.

The case of a Hilbert space in the domain was investigated in papers \cite{G-2024,G-2026}, 
where it is proved in \cite[Theorem 2.1]{G-2026} (see also \cite[Lemma 1.2]{G-2024}) that:

\begin{enumerate}[-]
\item\ 
an operator $T:{\cal H}\to Y$ from a Hilbert space ${\cal H}$ to an LCS $Y$ is topologically bounded if and only if $T$ 
takes orthonormal sequences onto topologically bounded ones.
\end{enumerate}

Let us drop requirements of additional structures (e.g., lattice operations, order, or inner products) on a space in the domain and consider:

\begin{question}\label{que 1}
{\em Is an operator $T:X\to Y$ between Banach spaces $X$ and $Y$ bounded, 
whenever it takes normalized (unconditional) basic sequences of $X$ onto bounded sequences of $Y$?}
\end{question}

\noindent
The ``conditional part'' of Question \ref{que 1} is open even when $X$ is reflexive, whereas the ``unconditional part'' 
has a negative answer due to the example of a reflexive Banach space with no unconditional basic sequence constructed by Gowers and Maurey \cite{GM-1993}
(every linear functional on this space is bounded on normalized unconditional basic sequences simply because of absence of such sequences,
so, we take any unbounded functional on this space). Therefore, it makes sense to deal with the following questions.

\begin{enumerate}[-]
\item\
Under which conditions on Banach spaces $X$ and $Y$ every operator $T:X\to Y$ that is bounded on normalized unconditional basic sequences of $X$ is bounded?
\item\
How to describe the class $S$ of Banach lattices such that every linear functional on a Banach lattice $X$ from $S$ is bounded, whenever it is bounded
on every normalized disjoint sequence of $X$? An example, recently constructed by Storozhuk in \cite{S-2025} suggests that many of Banach lattices 
(e.g., infinite dimensional discrete Banach function spaces) do not belong to the class $S$.
\end{enumerate}

\medskip
The present paper deals with operators between vector spaces, where the order structure is assumed in the domain or/and codomain. 
We focus on Banach spaces and, when suitable, Fr\'echet spaces in the domain. The case of Fr\'echet spaces in the domain 
occurred recently in study of compactly classified operators along disjoint bounded sequences \cite{G-2026}.

\medskip
This paper is organized as follows. In Section 2, we survey briefly some recently published results on order-to-topology continuous/boun\-ded operators 
between OBSs and TVSs. Then, we prove Theorem \ref{theorem 1} on continuity of order-to-topology bounded operators from an ordered Fr\'echet space 
with a closed generating cone to a TVS. In Section 3, we study conditions providing continuity of Levi and Lebesgue operators.
For instance, we prove: Theorem \ref{theorem 2} telling us that each quasi $\sigma$-Levi operator from an ordered Fr\'echet space with a closed generating cone 
to an ordered TVS with a normal cone is continuous; and Theorem \ref{theorem 5} saying that each $\sigma$-\text{\rm w}-Lebesgue operator
from an ordered Banach space with a closed generating normal cone to a normed space is bounded.

\medskip
For further unexplained terminology and notation, we refer the reader to \cite{ABo-2006,AB-2006,AT-2007}.

\section{On continuity of some operators between ordered topological vector spaces}

\hspace{4mm}
Order-to-topology continuous and order-to-norm bounded operators have been studied recently (mostly in the Banach lattice setting) by several authors 
(see, e.g., \cite{AEG-2022,E2-2026,EEG1-2025,G-2024,G-2026,JAM-2021,KTA-2025}). Some of their results were extended 
to operators from an OBS to a TVS in \cite{E1-2026}. It is important to find further conditions on operators in ordered TVSs which provide automatic 
continuity/boundedness. In this section, we continue the research initiated in \cite{E1-2026,E2-2026,E4-2025,EEG1-2025}. 
First, we point out Propositions \ref{prop 1} and \ref{prop 3} providing several conditions for topological boundedness 
of order-to-topology continuous and of order-to-norm bounded operators, and the Uniform Boundedness Principle 
for collectively qualified families of such operators. Then, we prove Theorem \ref{theorem 1} telling us that every 
order-to-topology bounded operator from an ordered Fr\'echet space with a closed generating cone to a TVS is continuous.
We need several further notions and notation. 

\begin{definition}\label{defi 1}
{\em A net $(x_\alpha)$ in an ordered vector space $X$ is
\begin{enumerate}[$a)$]
\item\ 
{\em order convergent} to $x$ (o-{\em convergent} to $x$, or $x_\alpha{\xrightarrow[]{o}} x$) if there exists a net
$g_\beta\downarrow 0$ in $X$ such that, for each $\beta$, there is $\alpha_\beta$ 
such that $\pm(x_\alpha-x)\le g_\beta$ for $\alpha\ge\alpha_\beta$. 
\item\ 
{\em relative uniform convergent} to $x$ (ru-{\em convergent} to $x$, or $x_\alpha{\xrightarrow[]{ru}}x$) 
if, for some $u\in X_+$, there exists an increasing sequence 
$(\alpha_n)$ of indices satisfying $\pm(x_\alpha-x)\le\frac{1}{n}u$ for $\alpha\ge\alpha_n$. 
\item\ 
{\em order bounded} if there exists an order interval $[a,b]$ in $X$ such that $x_\alpha\in[a,b]$ for all $\alpha$.
\item\ 
{\em eventually order bounded} if there exist $a,b\in X$ and $\alpha_0$ with $x_\alpha\in[a,b]$ for $\alpha\ge\alpha_0$.
\item\ 
o-{\em Cauchy} (resp., ru-{\em Cauchy}) if $(x_{\alpha'}-x_{\alpha''}){\xrightarrow[]{o}} 0$ (resp., $(x_{\alpha'}-x_{\alpha''}){\xrightarrow[]{ru}}0$). 
\end{enumerate}}
\end{definition}

\noindent
It is straightforward to see that every o-convergent and every ru-convergent (hence, every o-Cauchy and every ru-Cauchy) net is eventually order bounded.
In particular, o-convergent and ru-convergent (hence, also o-Cauchy and ru-Cauchy) sequences are order bounded.

\begin{definition}\label{defi 2}
{\em An operator $T$ from an OVS $X$ to an OVS $Y$ is
\begin{enumerate}[$a)$]
\item\
{\em order bounded} $(T\in{\cal L}_{ob}(X,Y))$ if $T[a,b]$ is order bounded for every $[a,b]\subseteq X$.
\item\
{\em order continuous} $(T\in{\cal L}_{oc}(X,Y))$ if $Tx_\alpha{\xrightarrow[]{o}} 0$ in $Y$, whenever $x_\alpha{\xrightarrow[]{o}} 0$.
\item\
{\em ru-continuous}  $(T\in{\cal L}_{rc}(X,Y))$ if $Tx_\alpha{\xrightarrow[]{ru}}0$ in $Y$, whenever $x_\alpha{\xrightarrow[]{ru}}0$.
\end{enumerate}
An operator $T$ from an OVS $X$ to a TVS $(Y,\tau)$ is
\begin{enumerate}
\item[$d)$]\
{\em order-to-topology continuous} $(T\in{\cal L}_{otc}(X,Y))$ if $x_\alpha{\xrightarrow[]{o}} 0\Longrightarrow Tx_\alpha{\xrightarrow[]{\tau}}0$.
\item[$e)$]\
{\em ru-to-topology continuous} $(T\in{\cal L}_{rtc}(X,Y))$ if $x_\alpha{\xrightarrow[]{ru}}0\Longrightarrow Tx_\alpha{\xrightarrow[]{\tau}}0$.
\item[$f)$]\
{\em order-to-topology bounded} $(T\in{\cal L}_{otb}(X,Y))$ if $T[a,b]$ is $\tau$-bounded for every $a,b\in X$.
\end{enumerate}
An operator $T$ from a TVS $(X,\xi)$ to a TVS $(Y,\tau)$ is
\begin{enumerate}[]
\item[$g)$]\ 
{\em continuous} $(T\in{\cal L}_{tc}(X,Y))$ if $Tx_\alpha{\xrightarrow[]{\tau}}0$, whenever $x_\alpha{\xrightarrow[]{\xi}}0$.
\item[$e)$]\ 
{\em topologically-bounded} $(T\in{\cal L}_{b}(X,Y))$ if $TU$ is $\tau$-bounded for every $\xi$-bounded $U$. 
\end{enumerate}}
\end{definition}

\noindent
Evidently, the sets ${\cal L}_{ob}(X,Y)$, ${\cal L}_{otc}(X,Y)$, etc. given by Definition \ref{defi 2} are vector spaces.

\medskip
We recall several notions on collectively qualified sets of operators 
(cf. \cite{E1-2025,E2-2025,E1-2026,E2-2026,G-2024,G-2026,EEG1-2025}). A family ${\cal B}=\{(x^b_\alpha)_{\alpha\in A}\}_{b\in B}$ 
of nets indexed by a directed set $A$ in a TVS $(X,\tau)$ {\em collective $\tau$-converges} to $0$ (shortly, ${\cal B}{\xrightarrow[]{c\tau}}0$) 
if, for each $U\in\tau(0)$, there exists $\alpha_U$ such that $x^b_\alpha\in U$ for all $\alpha\ge\alpha_U$ and $b\in B$. 

\begin{definition}\label{defi 3}
{\em A set ${\cal T}$ of operators from an OVS $X$ to an OVS $Y$ is 
\begin{enumerate}[$a)$]
\item\ 
{\em equi-order bounded} $({\cal T}\in\text{\bf L}_{ob}(X,Y))$, whenever ${\cal T}[a,b]$ is order bounded for every $a,b\in X$. 
\end{enumerate}
A set ${\cal T}$ of operators from an OVS $X$ to a TVS $(Y,\tau)$ is 
\begin{enumerate}
\item[$b)$]\
{\em equi-order-to-topology continuous} $({\cal T}\in\text{\bf L}_{otc}(X,Y))$ if 
$x_\alpha{\xrightarrow[]{o}} 0\Longrightarrow{\cal T}x_\alpha{\xrightarrow[]{c\tau}}0$.
\item[$c)$]\
{\em equi-\text{\rm ru}-to-topology continuous} $({\cal T}\in\text{\bf L}_{rtc}(X,Y))$ 
if $x_\alpha{\xrightarrow[]{ru}}0\Longrightarrow{\cal T}x_\alpha{\xrightarrow[]{c\tau}}0$.
\item[$d)$]\
{\em equi-order-to-topology bounded} $({\cal T}\in\text{\bf L}_{otb}(X,Y))$ if ${\cal T}[a,b]$ is $\tau$-bounded for every $a,b\in X$.
\end{enumerate}
A set ${\cal T}$ of operators from a TVS $(X,\xi)$ to a TVS $(Y,\tau)$ is 
\begin{enumerate}[]
\item[$e)$]\ 
{\em equi-continuous} $({\cal T}\in\text{\bf L}_{tc}(X,Y))$ if ${\cal T}x_\alpha{\xrightarrow[]{c\tau}}0$, whenever $x_\alpha{\xrightarrow[]{\xi}}0$. 
\item[$f)$]\ 
{\em equi-topologically-bounded} $({\cal T}\in\text{\bf L}_{b}(X,Y))$ if ${\cal T}U$ is $\tau$-bounded, whenever $U$ is $\xi$-bounded.
\end{enumerate}}
\end{definition}

\noindent
The sequential versions of Definition \ref{defi 2}\,b)--e),\,g) and Definition \ref{defi 3}\,b),\,c), and\,e) are obtained via replacing nets by sequences
(with notations ${\cal L}_{oc}^\sigma(X,Y)$, ${\cal L}_{otc}^\sigma(X,Y)$, $\text{\bf L}_{tc}^\sigma(X,Y)$, etc.).

\bigskip
We recall several recent results concerned the notions listed in Definitions \ref{defi 2} and \ref{defi 3}.

\begin{proposition}\label{prop 1}
{\em (\cite[Theorem 2.1]{E1-2026})}
If $X$ is an OVS and $Y$ is a TVS then $\text{\bf L}_{otb}(X,Y)\subseteq\text{\bf L}_{rtc}(X,Y)$. 
If additionally $X_+$ is generating then $\text{\bf L}_{otb}(X,Y)=\text{\bf L}_{rtc}(X,Y)$.
\end{proposition}

\noindent
In particular, each order-to-topology bounded operator from an OVS $X$ to a TVS $Y$ is ru-to-topology continuous,
and if additionally $X_+$ is generating, then ${\cal L}_{otb}(X,Y)={\cal L}_{rtc}(X,Y)$.
The following proposition is a partial extension of \cite[Theorem 2.1]{E2-2026} and \cite[Theorem 2.8]{EEG1-2025}.

\begin{proposition}\label{prop 3}
{\em (\cite[Theorem 2.5]{E1-2026})}
If $X$ is an OBS with a closed generating cone and $Y$ is a TVS then $\text{\bf L}_{otb}(X,Y)\subseteq\text{\bf L}_{b}(X,Y)$.
\end{proposition}

\noindent
It follows that every equi-order-to-norm bounded set of operators from an OBS with a closed generating cone to an NS is uniformly bounded.
In particular, each order bounded operator from an OBS with a closed generating cone to an ONS with a normal cone is continuous.
Following \cite{E1-2026}, we collect few more statements below.

\medskip
Let $X$ be an ordered Banach space with a closed generating cone.
\begin{enumerate}[$i)$]
\item\ 
If $Y$ is a Banach space, then an operator $T:X\to Y$ is continuous whenever it takes order intervals onto relatively weak compact sets \cite[Corollary 2.9]{E1-2026}.
\item\ 
If $X_+$ is normal and $Y$ is a NS, then $\text{\bf L}_{rtc}(X,Y)=\text{\bf L}_{otb}(X,Y)=\text{\bf L}_{b}(X,Y)$ \cite[Corollary 2.12]{E1-2026}.
\item\ 
If $X_+$ is normal and $Y$ is an ONS with a generating normal cone, then ${\cal L}_{rc}(X,Y)\subseteq{\cal L}_{b}(X,Y)$ \cite[Corollary 2.14]{E1-2026}.
\item\ 
If $Y$ is an OBS with a closed generating cone and $X^\sim\neq\emptyset$, where $X^\sim$ is an order dual of $X$,
then every ${\cal T}\in\text{\bf L}_{ob}(X,Y)$ is uniformly bounded if and only if $Y_+$ is normal \cite[Proposition 2.15]{E1-2026}.
\item\ 
If $X_+$ is normal and $Y=(Y,\tau)$ is a Banach space with the weak topology, 
then ${\cal L}_{otc}^\sigma(X,Y)\subseteq{\cal L}_{otb}(X,Y)$ \cite[Proposition 2.16]{E1-2026}.
\item\ 
If $X_+$ is closed generating and normal and $Y=(Y,\tau)$ is a Banach space with the weak topology, then
every order-to-topology $\sigma$-continuous operator $T:X\to Y$ is topologically bounded \cite[Corollary 2.18]{E1-2026}.
\end{enumerate}

\medskip
Now, we present several new results in the setting of Fr\'echet spaces in the domain.
Recall that a Fr\'echet space $({\cal F},\tau)$ is an LCS, whose topology $\tau$ is metrizable by a complete metric 
(we shall suppose that $\tau$ is induced by a metric $p(x,y)=\sum\limits_{k=1}^{\infty}\frac{2^{-k}p_k(x-y)}{1+p_k(x-y)}$ for $x,y\in Y$,
where $\{p_k\}_{k\in\mathbb{N}}$ is a separating family of semi-norms on ${\cal F}$). The following theorem generalizes \cite[Theorem 2.3]{E1-2026}.

\begin{theorem}\label{theorem 1}
$\text{\bf L}_{otb}({\cal F},Y)\subseteq\text{\bf L}_{tc}({\cal F},Y)$, 
whenever $({\cal F},\tau)$ is an ordered Fr\'echet space with a $\tau$-closed generating cone and $(Y,\xi)$ is a topological vector space.
\end{theorem}

\begin{proof}
Suppose, on the contrary, ${\cal T}\in\text{\bf L}_{otb}({\cal F},Y)\setminus\text{\bf L}_{tc}({\cal F},Y)$. 
For each $n\in\mathbb{N}$, denote 
$$
   V_n:=B_{\cal F}(0,2^{-2n})\cap{\cal F}_+-B_{\cal F}(0,2^{-2n})\cap{\cal F}_+.
$$
Clearly, $\{B_{\cal F}(0,2^{-2n}):n\in\mathbb{N}\}$ is a $\tau$-neighborhood base of zero consisting of $\tau$-closed sets 
and $B_{\cal F}(0,2^{-2(n+1)})+B_{\cal F}(0,2^{-2(n+1)})\subseteq B_{\cal F}(0,2^{-2n})$ holds for each $n$.
Furthermore, each $B_{\cal F}(0,2^{-2n})$ is circled. To see this, it suffices to show that $p(x,0)\le 2^{-2n}$ implies $p(tx,0)\le 2^{-2n}$ whenever $|t|\le 1$.  
So, let $p(x,0)\le 2^{-2n}$ and $|t|\le 1$. Since $\frac{|t|a}{1+|t|a}\le\frac{a}{1+a}$ for every $a\in\mathbb{R}_+$, we obtain
$p(tx,0)=\sum\limits_{k=1}^{\infty}\frac{2^{-k}p_k(tx)}{1+p_k(tx)}=\sum\limits_{k=1}^{\infty}\frac{2^{-k}|t|p_k(x)}{1+|t|p_k(x)}\le
\sum\limits_{k=1}^{\infty}\frac{2^{-k}p_k(x)}{1+p_k(x)}=p(x,0)\le 2^{-2n}$.
By the Ando theorem \cite[Theorem 2.10]{AT-2007}, $\{V_n:n\in\mathbb{N}\}$ is a $\tau$-neighborhood base of zero
consisting of circled sets satisfying $V_{n+1}+V_{n+1}\subseteq V_n$ for all $n$.

\medskip
Since ${\cal T}\notin\text{\bf L}_{tc}({\cal F},Y)$, it follows from \cite[Theorem 5.6]{ABo-2006} that there exist a neighborhood base $\cal{B}$ of $\xi$ at zero
consisting of circled sets, a sequence $(y_n)$ with $y_n\in V_n$, and a sequence $(T_n)$ in ${\cal T}$ such that $T_ny_n\notin W+W$ for 
for some $W\in\cal{B}$. Find $y'_n$ and $y''_n$ in $B_{\cal F}(0,2^{-2n})\cap{\cal F}_+$ with $y_n=y'_n-y''_n$ for all $n$.
Since $T_ny'_n-T_ny''_n=T_ny_n\notin W+W$, pick $z_n=y'_n$ if $T_ny'_n\notin W$ or $z_n=y''_n$ otherwise (then, $T_ny''_n\notin W$). 

\medskip
Now, let $x_n=2^nz_n$ for each $n\in\mathbb{N}$. It follows from $p(z_n,0)\le 2^{-2n}$ that $p(x_n,0)\le 2^{-n}$ for all $n\in\mathbb{N}$.
Take $x\in{\cal F}$ so that $p\bigl(x-\sum\limits_{n=1}^k x_n\bigl)\to 0$ as $k\to\infty$.
Such $x$ exists because $\big(\sum\limits_{n=1}^k x_n\big)_k$ is a $p$-Cauchy sequence and the metric $p$ is complete.
As $x_n\in{\cal F}_+$ for all $n$, the $p$-closeness of ${\cal F}_+$ ensures $x\ge 0$ and $x_n\in[0,x]$ for all $n$.
Since ${\cal T}\in\text{\bf L}_{otb}({\cal F},Y)$ then ${\cal T}[0,x]\subseteq NW$ for some $N\in\mathbb{N}$.
It follows from $x_n\in[0,x]$ that $T_nx_n\in NW$ for every $n$.
This is absurd, because $T_nx_n=T_n(2^nz_n)=2^nT_nz_n\notin 2^nW$ for all $n$.
\end{proof}

\begin{corollary}\label{coro 1}
Every order-to-topology bounded operator from an ordered Fr\'echet space with a closed generating cone to a topological vector space is continuous.
\end{corollary}

\noindent
Keeping in mind that ${\cal L}_{ob}(X,Y)\subseteq{\cal L}_{otb}(X,Y)$ if $Y$ is a TVS with a normal cone \cite[Lemma 2.22]{AT-2007},
we obtain the following conditions for automatic continuity of order bounded operators.

\begin{corollary}\label{coro 2}
Every order bounded operator from an ordered Fr\'echet space with a closed generating cone to an 
ordered topological vector space with a normal cone is continuous.
\end{corollary}

\section{On continuity of the Levi and the Lebesgue operators}

\hspace{4mm}
In this section, we study $\tau$-Levi and $\tau$-Lebesgue operators \cite{AEG-2022,GE-2022,E1-2026} in the more general OVS-setting. 
We prove Theorem \ref{theorem 2} saying that a quasi $\sigma$-Levi operator from an ordered Fr\'echet space with a closed generating 
cone to an ordered topological vector space with a normal cone is continuous. This theorem is a far going generalization of \cite[Theorem 3.1]{E4-2025},
where the continuity is established for a quasi $\sigma$-Levi operator from a Banach lattice to a normed lattice.
We mention also: Lemma \ref{lemma 3} saying that each $\sigma$-\text{\rm w}-Lebesgue operator is order-to-norm bounded whenever an OBS in the domain has a normal cone;
Theorem \ref{theorem 5} saying that, if the domain has closed generating normal cone, then every $\sigma$-\text{\rm w}-Lebesgue operator is bounded; and
Theorem \ref{theorem 3} providing conditions for ($\sigma$-)order to (weak) norm continuity of ($\sigma$-)(\text{\rm w}-)Lebesgue operators.
We shall use below the following versions of $\tau$-Levi and $\tau$-Lebesgue operators 
(see \cite{AEG-2022,E4-2025,GE-2022,JAM-2021} in the vector lattice-setting).

\begin{definition}\label{defi 4}
{\em An operator $T$ from an ordered TVS $(X,\tau)$ to an OVS $Y$ is
\begin{enumerate}[$a)$]
\item\
{\em quasi $\sigma$-Levi} if $T$ takes $\tau$-bounded increasing sequences of $X$ to o-Cau\-chy ones.
\end{enumerate}
An operator $T$ from an OVS $X$ to an LCS $(Y,\xi)$ is 
\begin{enumerate}
\item[$b)$]\
$(\sigma$-$)${\em Lebesgue} if $Tx_\alpha{\xrightarrow[]{\tau}}0$ for every $($sequence$)$ net $x_\alpha\downarrow 0$ in $X$.
The set of such operators is denoted by ${\cal L}_{Leb}(X,Y)$ (resp., ${\cal L}^\sigma_{Leb}(X,Y)$). 
\item[$c)$]\
$(\sigma$-$)$\text{\rm w}-{\em Lebesgue} if $Tx_\alpha{\xrightarrow[]{\sigma(Y,Y^\prime)}} 0$ for every $($sequence$)$ net $x_\alpha\downarrow 0$ in $X$.
The set of such operators is denoted by ${\cal L}_{wLeb}(X,Y)$ (resp., ${\cal L}^\sigma_{wLeb}(X,Y)$). 
\item[$d)$]\
$(\sigma$-$)${\em order-to-weak continuous} if $Tx_\alpha{\xrightarrow[]{\sigma(Y,Y^\prime)}} 0$ for every $($sequence$)$ net $x_\alpha{\xrightarrow[]{o}} 0$ in $X$.
The set of such operators is denoted by ${\cal L}_{owc}(X,Y)$ (resp., ${\cal L}^\sigma_{owc}(X,Y)$) (cf., Definition \ref{defi 2}\,$d)$). 
\end{enumerate}}
\end{definition}

\medskip
First, we prove that under rather mild conditions, a quasi $\sigma$-Levi operator is continuous.
This generalizes \cite[Theorem 3.1]{E4-2025} which states that every quasi $\sigma$-Levi operator from a Banach lattice to a normed lattice is continuous.

\begin{theorem}\label{theorem 2}
Every quasi $\sigma$-Levi operator from an ordered Fr\'echet space with a closed generating cone to an ordered topological vector space with a normal cone is continuous.
\end{theorem}

\begin{proof}\
Let $T:{\cal F}\to Y$ be a quasi $\sigma$-Levi operator from an ordered Fr\'echet space $({\cal F},\tau)$ with a $\tau$-closed generating cone
to an ONS $Y$ with a normal cone. Toward a contradiction, suppose that $T$ is not continuous.

\medskip
Arguing as in the proof of Theorem \ref{theorem 1},
consider a $\tau$-neighborhood base $\{B_{\cal F}(0,2^{-2n})\}_{n\in\mathbb{N}}$ of zero consisting of $\tau$-closed circled sets,
and denote $V_n:=B_{\cal F}(0,2^{-2n})\cap{\cal F}_+-B_{\cal F}(0,2^{-2n})\cap{\cal F}_+$ for each $n\in\mathbb{N}$.
Then $\{V_n:n\in\mathbb{N}\}$ is a $\tau$-neighborhood base of zero consisting of circled subsets satisfying $V_{n+1}+V_{n+1}\subseteq V_n$ for all $n$.
Since $T$ is not continuous, there exist a neighborhood base $\cal{B}$ of $\xi$ at zero
consisting of circled sets, a sequence $(y_n)$ with $y_n\in V_n$, and a sequence $(T_n)$ in ${\cal T}$ such that $T_ny_n\notin W+W$ for for some $W\in\cal{B}$. 
As in the proof of Theorem \ref{theorem 1}, for each $n\in\mathbb{N}$ choose $z_n\in B_{\cal F}(0,2^{-2n})\cap{\cal F}_+$ with $T_nz_n\notin W$.

\medskip
Letting $x_n:=2^nz_n$, denote $u_i=\sum\limits_{k=1}^i x_k$ for each $i\in\mathbb{N}$. Clearly, $(u_i)$ is an increasing sequence. 
Take an $l\in\mathbb{N}$. It is straightforward to see that, for each $m\in\mathbb{N}$, 
$$
   \sum\limits_{k=l+1}^{m+l+1}x_k=\sum\limits_{k=l+1}^{m+l+1}2^kz_k\in 
   \sum\limits_{k=l+1}^{m+l+1}2^kB_{\cal F}(0,2^{-2k})\subseteq 
   \sum\limits_{k=0}^{m}2^{-k}B_{\cal F}(0,2^{-(l+1)})\subseteq B_{\cal F}(0,2^{-l}).
$$
Find an $N>0$ satisfying $\sum\limits_{k=1}^{j}x_k\in NB_{\cal F}(0,2^{-l})$ for $j=1,...,l$. Then, for every $i\in\mathbb{N}$, 
$$
   u_{i}=\sum\limits_{k=1}^{i}x_k=\sum\limits_{k=1}^{\min(i,l)}x_k+\sum\limits_{k=\min(i,l)+1}^{i}x_k\in (N+1)B_{\cal F}(0,2^{-l}).
$$ 
As $\{B_{\cal F}(0,2^{-l})\}_{l\in\mathbb{N}}$ is a $\tau$-neighborhood base of zero, the sequence $(u_i)$ is $\tau$-bounded.

\medskip
As $T$ is a quasi $\sigma$-Levi operator, the sequence $(Tu_i)$ is o-Cauchy in $Y$, and hence is order bounded. 
Since $Y_+$ is normal, the sequence $(Tu_i)$ is topologically bounded by \cite[Lemma 2.22]{AT-2007}. 
It is absurd because $Tu_{i+1}-Tu_i=Tx_{i+1}=2^{i+1}Tz_{i+1}\notin 2^{i+1}W$ for all $i\in\mathbb{N}$.
\end{proof}

\medskip
Recall that a positive operator $T$ from a vector lattice to a normed lattice is Lebesgue/$\sigma$-Lebesgue 
if and only if $T$ is order-to-norm continuous/$\sigma$-order-to-norm continuous \cite[Lemma 2.1]{AEG-2022}. 
The same result holds also in the following more general setting.

\begin{theorem}\label{theorem 3}
Let $T$ be a positive operator from an ordered vector space $X$ to an ordered topological vector space $(Y,\tau)$ with a normal cone. Then
$$
   T\in{\cal L}_{Leb}(X,Y)\Longleftrightarrow T\in{\cal L}_{otc}(X,Y) \ \ \ and \ \ \
   T\in{\cal L}^\sigma_{Leb}(X,Y)\Longleftrightarrow T\in{\cal L}^\sigma_{otc}(X,Y).
$$
\end{theorem}

\begin{proof}
As the $\sigma$-Lebesgue case is similar, we consider only the first equivalence. Since obviously ${\cal L}_{otc}(X,Y)\subseteq{\cal L}_{Leb}(X,Y)$,
it suffices to show that $T\in{\cal L}_{Leb}(X,Y)$ implies $T\in{\cal L}_{otc}(X,Y)$.

On the way to a contradiction, suppose $T\in{\cal L}_{Leb}(X,Y)\setminus{\cal L}_{otc}(X,Y)$.
Then, there exist a net $x_\alpha{\xrightarrow[]{o}} 0$ in $X$ and a neighborhood $U\in\tau(0)$ such that, for every $\alpha$, there exists $\alpha'\ge\alpha$
with $Tx_{\alpha'}\notin U$. Since $Y_+$ is normal, we may suppose that $U$ is circled and full by \cite[Corollary 2.21]{AT-2007}. 

Pick a net $g_\beta\downarrow 0$ in $X$ such that, for each $\beta$, there exists $\alpha_\beta$
with $x_\alpha\in[-g_\beta,g_\beta]$ for $\alpha\ge\alpha_\beta$.  
As $T$ is Lebesgue, $Tg_\beta{\xrightarrow[]{\tau}}0$. Since $T\ge 0$, $Tg_\beta\downarrow\ge 0$.
Suppose $Tg_\beta\ge f\ge 0$ for some $f\in Y$.
Since $\tau$ is a Hausdorff topology, it follows from \cite[Theorem 2.23]{AT-2007} that $f=0$, and hence $Tg_\beta\downarrow 0$.
Since $T\ge 0$, we obtain $Tx_{\alpha}\in[-Tg_\beta,Tg_\beta]$ for $\alpha\ge\alpha_\beta$. 
Since $Tg_\beta{\xrightarrow[]{\tau}}0$, there exists $\beta_U$ such that $Tg_\beta\in U$ for all $\beta\ge\beta_U$.
Since $U$ is circled and full, $[-Tg_{\beta_U},Tg_{\beta_U}]\subseteq U$ for all $\beta\ge\beta_U$, and hence $Tx_{\alpha}\in U$  
for all $\alpha\ge\alpha_{\beta_U}$. It is absurd because $Tx_{\alpha'_{\beta_U}}\notin U$.
Thus, $T\in{\cal L}_{otc}(X,Y)$.
\end{proof}

\begin{corollary}\label{coro 3}
Let $X$ be an OVS and $Y$ be a normal ONS. Then $({\cal L}_{Leb})_+(X,Y)=({\cal L}_{otc})_+(X,Y)$ and 
$({\cal L}^\sigma_{Leb})_+(X,Y)=({\cal L}^\sigma_{otc})_+(X,Y)$.
\end{corollary}

\medskip
By \cite[Proposition 2.2]{AEG-2022}, a positive order continuous operator $T$ from a vector lattice to a normed lattice 
is Lebesgue if and only if $T$ is \text{\rm w}-Lebesgue. We extend this fact to the OVS-setting as follows.

\begin{theorem}\label{theorem 4}
Let $T$ be a positive order continuous operator from an OVS $X$ to a Hausdorff LCS $(Y,\tau)$ ordered by a normal cone. Then 
$$
   T\in{\cal L}_{Leb}(X,Y)\Longleftrightarrow T\in{\cal L}_{wLeb}(X,Y) \ \ \ and \ \ \
   T\in{\cal L}^\sigma_{Leb}(X,Y)\Longleftrightarrow T\in{\cal L}^\sigma_{wLeb}(X,Y).
$$
\end{theorem}

\begin{proof}
As the $\sigma$-Lebesgue case is similar, we consider only the first equivalence. 
Suppose $T\in{\cal L}_{wLeb}(X,Y)$ and take $x_\alpha\downarrow 0$ in $X$.
Since $T\ge 0$, we have $Tx_\alpha\downarrow\ge 0$ and, as $T$ is \text{\rm w}-Lebesgue, $Tx_\alpha{\xrightarrow[]{\sigma(Y,Y^\prime)}} 0$. 
By \cite[Lemma 2.28]{AT-2007}, we conclude $Tx_\alpha{\xrightarrow[]{\tau}}0$. Thus, $T\in{\cal L}_{Leb}(X,Y)$.
\end{proof}

\begin{corollary}\label{coro 4}
Let $X$ be an OVS and $Y$ be a normal ONS. Then $({\cal L}_{Leb})_+(X,Y)=({\cal L}_{wLeb})_+(X,Y)$ and 
$({\cal L}^\sigma_{Leb})_+(X,Y)=({\cal L}^\sigma_{wLeb})_+(X,Y)$.
\end{corollary}

\begin{question}\label{que 2}
{\em Are the conclusions of Corollaries \ref{coro 3} and \ref{coro 4} hold true for arbitrary (not necessarily positive) operators?}
\end{question}

\medskip
The following lemma is a direct consequence of Definition \ref{defi 4}.

\begin{lemma}\label{lemma 1}
Let $X$ be an OVS and $Y$ be an LCS. Then 
$$
   {\cal L}_{Leb}(X,Y)\subseteq{\cal L}_{wLeb}(X,Y), \ \ 
   {\cal L}^\sigma_{Leb}(X,Y)\subseteq{\cal L}^\sigma_{wLeb}(X,Y), 
$$
$$
   {\cal L}_{otc}(X,Y)\subseteq{\cal L}_{owc}(X,Y), \ \ 
   {\cal L}^\sigma_{otc}(X,Y)\subseteq{\cal L}^\sigma_{owc}(X,Y),  
$$
$$
   {\cal L}_{otc}(X,Y)\subseteq{\cal L}_{Leb}(X,Y)\bigcap{\cal L}^\sigma_{otc}(X,Y)\subseteq 
   {\cal L}_{Leb}(X,Y)\bigcup{\cal L}^\sigma_{otc}(X,Y)\subseteq{\cal L}^\sigma_{Leb}(X,Y), \ \ \text{\rm and}
$$
$$
   {\cal L}_{owc}(X,Y)\subseteq{\cal L}_{wLeb}(X,Y)\bigcap{\cal L}^\sigma_{owc}(X,Y)\subseteq 
   {\cal L}_{wLeb}(X,Y)\bigcup{\cal L}^\sigma_{owc}(X,Y)\subseteq{\cal L}^\sigma_{wLeb}(X,Y).
$$
\end{lemma}

\begin{lemma}\label{lemma 2}
Let $X$ be an OBS with a closed generating normal cone and let $Y$ be an NS. Then, ${\cal L}_{otb}(X,Y)={\cal L}_{b}(X,Y)$.
\end{lemma}

\begin{proof}
The inclusion ${\cal L}_{otb}(X,Y)\subseteq{\cal L}_{b}(X,Y)$ follows from \cite[Corollary 2.9]{EEG1-2025}. The reverse inclusion hods due to normality of $X$.
\end{proof}

Another case of automatic continuity of operators concerns their order continuity. 
Namely, as every order continuous operator from an Archimedean OVS with generating cone to an OVS is order bounded \cite[Corollary 2.2]{EEG1-2025}, we have:
${\cal L}_{oc}(X,Y)\subseteq{\cal L}_{b}(X,Y)$, whenever $X$ and $Y$ are OBSs with generating cones such that $X_+$ is closed and $Y_+$ is normal.
Condition $X_+-X_+=X$ is essential here since each $T:X\to X$ is order bounded and order continuous whenever $X_+=\{0\}$.

\begin{example}\label{example 1}
{\em (\cite[Example 2.12 b)]{EEG1-2025})
Choose a normalized Hamel basis $\{\xi:\xi\in\Xi\}$ of an infinite dimensional Banach space $X$ and identify $X$ with $c_{00}(\Xi)$.
Consider a positive cone $X_+:=\Big\{\sum\limits_{\xi\in\Xi}\lambda_\xi\xi: \lambda_\xi\ge 0\Big\}$ in $X$.
The cone $X_+$ is generating. However $X_+$ is not closed $($e.g., by \text{\rm \cite[Proposition 4.1]{AN-2009}}$)$. 
Pick any sequence $(\xi_n)_{n=0}^\infty$ of distinct elements of $\Xi$ and define 
$T:(X,X_+)\to X$ by $T\xi=0$ for $\xi\in\Xi\setminus\{\xi_n:n\in\mathbb{N}\}$ and $T\xi_n=n\xi_0$ for $n\in\mathbb{N}$.
It is straightforward to see that the operator $T$ is order-to-norm bounded and positive, but $T$ is not bounded.}
\end{example}

It is worth nothing that the order-to-norm and order-to-weak continuous operators on vector lattices were defined in \cite{JAM-2021}
without specifying the class of (just linear or continuous linear) operators. Automatic boundedness of such operators was investigated in the recent paper \cite{KTA-2025}. 

\medskip
We are going to establish further conditions for automatic boundedness of order-to-weak continuous operators between ordered Banach spaces and and normed spaces. 
More precisely, we prove Theorem \ref{theorem 5} which tells us that every $\sigma$-\text{\rm w}-Lebesgue operator from an OBS with 
a closed generating normal cone to an NS is bounded. Let us begin with conditions under which $\sigma$-\text{\rm w}-Lebesgue operators are order-to-norm bounded. 

\begin{lemma}\label{lemma 3}
Let $X$ be an ordered Banach space with a normal cone and $Y$ a normed space. Then ${\cal L}^\sigma_{wLeb}(X,Y)\subseteq{\cal L}_{otb}(X,Y)$. 
In particular, ${\cal L}^\sigma_{owc}(X,Y)\subseteq{\cal L}_{otb}(X,Y)$.
\end{lemma}

\begin{proof}
Let $T:X\to Y$ be $\sigma$-\text{\rm w}-Lebesgue. If $T$ is not order-to-norm bounded then $T[0,u]$ is not bounded for some $u\in X_+$. 
Since $X$ is normal then $[0,u]$ is bounded, say $\sup\limits_{x\in[0,u]}\|x\|\le M$. 
Pick a sequence $(u_n)$ in $[0,u]$ with $\|Tu_n\|\ge n2^n$, and set $y_n:=\|\cdot\|$-$\sum\limits_{k=n}^\infty 2^{-k}u_k$ for $n\in\mathbb{N}$.
Then $y_n\downarrow\ge 0$. Let $0\le y_0\le y_n$ for all $n\in\mathbb{N}$.
Since $0\le y_0\le 2^{1-n}u$ it follows that $\|y_0\|\le M2^{1-n}$ for every $n$.
Thus, $y_n\downarrow 0$. It follows from $T\in{\cal L}^\sigma_{wLeb}(X,Y)$ that $(Ty_n)$ is \text{\rm w}-null, 
and hence is norm bounded, that is absurd because
$\|Ty_{n+1}-Ty_n\|=\|T(2^{-n}u_n)\|\ge n$ for all $n\in\mathbb{N}$. 
\end{proof}

\noindent
Combining Lemmas \ref{lemma 1}, \ref{lemma 2}, and \ref{lemma 3}, we obtain.

\begin{theorem}\label{theorem 5}
Let $X$ be an ordered Banach space with a closed generating normal cone and $Y$ a normed space. 
Then ${\cal L}^\sigma_{wLeb}(X,Y)\subseteq{\cal L}_{b}(X,Y)$. In particular, ${\cal L}^\sigma_{ows}(X,Y)\subseteq{\cal L}_{b}(X,Y)$.
\end{theorem}

\begin{corollary}\label{coro 5}
{\rm (cf. \cite[Proposition 3.2]{AEG-2022}, \cite[Lemma 2.3]{KTA-2025})}
Each $\sigma$-\text{\rm w}-Lebesgue, and hence each $\sigma$-order-to-weak continuous, operator from a Banach lattice to an normed space is bounded.
\end{corollary}

\noindent
Positive order-to-norm continuous order-to-norm bounded operators between normed lattices are not necessary bounded 
\cite[Remark 1.1.\,c)]{AEG-2022}. From the other hand, even a positive $\sigma$-Lebesgue rank one contraction 
on an \text{\rm AM}-space need not to be \text{\rm w}-Lebesgue \cite[Example~3.1]{AEG-2022}.

\begin{definition}\label{defi 5}
{\em The norm in an ordered normed space $X$ is 
{\em order continuous} if $\|x_\alpha\|\to 0$ for every net $x_\alpha{\xrightarrow[]{o}} 0$ in $X$.}
\end{definition}

\noindent
In certain cases o-convergence is equivalent to ru-convergence.

\begin{lemma}\label{lemma 4}
Let $X$ be an OBS with a closed normal cone and the norm in $X$ is order continuous. Then,
$x_\alpha{\xrightarrow[]{o}} 0\Longleftrightarrow x_\alpha{\xrightarrow[]{ru}}0$ in $X$.
\end{lemma}

\begin{proof}
Since $X_+$ is normal, $X$ is Archimedean, and hence $x_\alpha{\xrightarrow[]{ru}}0\Longrightarrow x_\alpha{\xrightarrow[]{o}} 0$ in $X$. 

\medskip
Let $X$ has order continuous norm and $x_\alpha{\xrightarrow[]{o}} 0$ in $X$. Take a net $g_\beta\downarrow 0$ in $X$
such that, for each $\beta$ there is  $\alpha_\beta$ satisfying $\pm x_\alpha\le g_\beta$ for $\alpha\ge\alpha_\beta$.
Since the norm in $X$ is order continuous and $X_+$ is normal, there exists
an increasing sequence $(\beta_k)$ such that $\sup\limits_{z\in[-g_{\beta_k},g_{\beta_k}]}\|z\|\le k^{-3}$ for all $k\in\mathbb{N}$.
Set $u=\|\cdot\|$-$\sum_{k=1}^\infty kg_{\beta_k}$. Then $\pm x_\alpha\le g_{\beta_k}\le k^{-1}u$ for $\alpha\ge\alpha_{\beta_k}$.
Choose any increasing sequence $(\alpha_k)$ with $\alpha_k\ge\alpha_{\beta_1},\alpha_{\beta_2},\ldots,\alpha_{\beta_k}$.
Then $\pm x_\alpha\le k^{-1}u$ for $\alpha\ge\alpha_k$, and hence $x_\alpha{\xrightarrow[]{ru}}0$.
\end{proof}

\begin{theorem}\label{theorem 6}
Let $X$ be an ordered Banach space with a closed generating normal cone and $Y$ a normed space.  
If $X$ has order continuous norm then 
$$
   {\cal L}_{Leb}(X,Y)={\cal L}_{otc}(X,Y) \ \ and \ \ {\cal L}_{wLeb}(X,Y)={\cal L}_{owc}(X,Y).
$$
\end{theorem}

\begin{proof}
In view of Lemma \ref{lemma 1}, we need to prove ${\cal L}_{Leb}(X,Y)\subseteq{\cal L}_{otc}(X,Y)$ 
and ${\cal L}_{wLeb}(X,Y)\subseteq{\cal L}_{owc}(X,Y)$. 
Let $x_\alpha{\xrightarrow[]{o}} 0$ in $X$. Then $x_\alpha{\xrightarrow[]{ru}}0$ by Lemma \ref{lemma 4}, and hence $\|x_\alpha\|\to 0$ due to normality of $X_+$. 

If $T\in{\cal L}_{Leb}(X,Y)$ then $T\in{\cal L}^\sigma_{Leb}(X,Y)$. By Theorem \ref{theorem 5},
$T\in{\cal L}_{b}(X,Y)$, and hence $\|Tx_\alpha\|\to 0$. Thus, $T\in{\cal L}_{otc}(X,Y)$. If $T\in{\cal L}_{wLeb}(X,Y)$ then 
again $T\in{\cal L}_{b}(X,Y)$ by Theorem \ref{theorem 5}. Thus, $\|Tx_\alpha\|\to 0$, and hence $Tx_\alpha{\xrightarrow[]{w}}0$.
So, $T\in{\cal L}_{owc}(X,Y)$.
\end{proof}

\medskip 
The next corollary is a partial extension of results of \cite[Lemma 1]{JAM-2021}, \cite[Lemma 2.1]{AEG-2022},  
\cite[Theorem 3]{ZSC-2023}, and \cite[Proposition 1.5]{KTA-2025} to the OVS setting.

\begin{corollary}\label{coro 7}
Let $E$ be a Banach lattice with order continuous norm and $Y$ a normed space.  Then,
${\cal L}_{Leb}(E,Y)={\cal L}_{otc}(E,Y)$ and ${\cal L}_{wLeb}(E,Y)={\cal L}_{owc}(E,Y)$.
\end{corollary}

\medskip
Finally, we give the conditions for order-to-topology continuity of order bounded operators.

\begin{theorem}\label{theorem 7}
Let $X$ be an ordered Banach space with a closed normal cone and order continuous norm, $Y$ be an ordered normed space with a generating normal cone,
and $T:X\to Y$ and order bounded operator. Then $T$ is order-to-topology continuous.
\end{theorem}

\begin{proof}
Let $x_\alpha{\xrightarrow[]{o}} 0$ in $X$. By Lemma \ref{lemma 4}, $x_\alpha{\xrightarrow[]{ru}}0$.
It follows from \cite[Corollary 2.5]{EEG1-2025} that $Tx_\alpha{\xrightarrow[]{ru}}0$. The normality of the cone $Y_+$ ensures $\|Tx_\alpha\|\to 0$.
\end{proof}

\bigskip
{}
\end{document}